\documentclass[11pt,a4paper]{amsart}
\usepackage[T1]{fontenc}
\usepackage{lmodern}
\usepackage{amsmath,amssymb,amsthm,mathtools}
\usepackage[margin=27mm]{geometry}
\usepackage{microtype}
\usepackage{enumitem}
\usepackage{xcolor}

\newif\ifshowchanges
\showchangestrue

\usepackage[hidelinks]{hyperref}
\newtheorem{theorem}{Theorem}[section]
\newtheorem{lemma}[theorem]{Lemma}
\newtheorem{proposition}[theorem]{Proposition}
\newtheorem{corollary}[theorem]{Corollary}
\newtheorem{conjecture}[theorem]{Conjecture}
\theoremstyle{definition}
\newtheorem{definition}[theorem]{Definition}
\theoremstyle{remark}
\newtheorem{remark}[theorem]{Remark}
\numberwithin{equation}{section}
\setlist[enumerate]{label=\textup{(\roman*)},leftmargin=*,itemsep=3pt,topsep=5pt}
\newcommand{\cT}{\mathcal T}
\newcommand{\cS}{\mathcal S}
\newcommand{\cW}{\mathcal W}
\DeclareMathOperator{\cl}{cl}
\DeclareMathOperator{\br}{br}
\newcommand{\doi}[1]{\href{https://doi.org/#1}{\nolinkurl{doi:#1}}}

\title{Connectivity of Matroids with Breadth-Critical Tangles}
\author[Zhuo Li]{Zhuo Li}
\address{School of Mathematical Sciences, Xiamen University}
\email{lzhuo@stu.xmu.edu.cn}
\author[Zhen Chen]{Zhen Chen$^{*}$}
\address{School of Mathematical Sciences, Xiamen University}
\email{chenzzhen@126.com}
\thanks{$^{*}$Corresponding author.}
\date{}

\begin{document}
\begin{abstract}
We construct an explicit infinite sequence of nonnegative integers
$(s_0,s_1,\ldots)$ such that, for every $k\ge2$, every
finite matroid with a breadth-critical tangle of order
at least $k$ is $(s_0,s_1,\ldots,s_{k-2})$-connected.
The constants $s_i$ depend only on $i$. This proves a conjecture of Brettell, Jowett, Oxley,
Semple, and Whittle.
\end{abstract}
\maketitle

\section{Introduction}
\label{sec:introduction}
Robertson and Seymour~\cite{RS91} introduced tangles to describe
highly connected regions of a graph and observed that the definition
extends to matroids. The matroid theory of tangles was
subsequently developed by Dharmatilake~\cite{Dha96}. A matroid with a tangle of large order need not
itself be highly connected. This motivates the study
of minors that retain a given tangle while satisfying
stronger connectivity conditions.

Brettell, Jowett, Oxley, Semple, and
Whittle~\cite{BJOSW25} introduced the \emph{breadth}
of a tangle $\cT$, defined as the maximum size of a
spanning uniform restriction of its associated
tangle matroid $M_{\cT}$.
A tangle $\cT$ of order $q$ in $M$ \emph{generates}
a tangle $\cS$ of order $q$ in a minor $N$ if
$\cS$ is the unique tangle of order $q$ in $N$
containing $A\cap E(N)$ for every $A\in\cT$.
Their main theorem states that every tangle of order
at least four generates a tangle of the same order
and breadth in a weakly $4$-connected minor.

We consider the connectivity of matroids with
breadth-critical tangles. A tangle $\cT$ of order
$q$ in $M$ is \emph{breadth-critical} if every
proper minor $N$ in which $\cT$ generates a
tangle $\cS$ of order $q$ satisfies
$\br_N(\cS)<\br_M(\cT)$. Brettell et al.~\cite[Corollaries~50, 53, and~66]{BJOSW25}
proved that a matroid with a breadth-critical tangle of order $q$
is connected if $q\ge2$, $3$-connected if $q\ge3$, and weakly
$4$-connected if $q\ge4$. What can be said at higher orders?
These results motivate a generalization of the
connectivity conditions above. We use the following
notation to formulate it.

For a matroid
$M$ with ground set $E(M)$ and rank function $r_M$, put
\[
 \lambda_M(A)=r_M(A)+r_M(E(M)\setminus A)-r(M)
 \qquad(A\subseteq E(M)),
\]
where $r(M)=r_M(E(M))$. For nonnegative integers $s_0,\ldots,s_t$,
a matroid $M$ is \emph{$(s_0,\ldots,s_t)$-connected} if, for every
$A\subseteq E(M)$ and every $i\in\{0,\ldots,t\}$,
\[
 \lambda_M(A)=i
 \quad\Longrightarrow\quad
 \min\{|A|,|E(M)\setminus A|\}\le s_i.
\]
In particular, weak $4$-connectivity is equivalent to
$(0,1,4)$-connectivity. 

Brettell et al.~\cite{BJOSW25} proposed the following conjecture.

\begin{conjecture}\label{conj:81}
There is an infinite sequence $(s_0,s_1,\ldots)$ of nonnegative
integers such that, for every $k\ge2$, if a finite matroid $M$
has a breadth-critical tangle of order at least $k$, then $M$ is
$(s_0,\ldots,s_{k-2})$-connected.
\end{conjecture}

Our main theorem proves the conjecture with an explicit sequence.

\begin{theorem}\label{thm:main}
For every integer $k\ge2$, every finite matroid
with a breadth-critical tangle of order at least $k$
is $(s_0,s_1,\ldots,s_{k-2})$-connected, where
$(s_i)_{i\ge0}$ is the sequence defined
in~\eqref{eq:recurrence} below.
\end{theorem}

The sequence is defined as follows. For integers $t\ge0$, let
\begin{equation}\label{eq:constants}
\begin{aligned}
 c_t&=\frac{3^t-1}{2},
 & b_t&=(2t+1)2^{2t+1},\\
 h_t&=c_t+\binom{2^{c_t}}{2}(b_t-1).
\end{aligned}
\end{equation}
Set $a_0=0$, and define
\begin{equation}\label{eq:recurrence}
\begin{aligned}
 a_t&=\max\{h_t,3a_{t-1}\} &&(t\ge1),\\
 s_i&=i+a_{2i} &&(i\ge0).
\end{aligned}
\end{equation}

The proof is based on a bound for the size of a $\cT$-weak set disjoint from a witness for the breadth of a breadth-critical tangle $\cT$. Using the linking theorem of Chen and Whittle~\cite{CW14}, we show that titanic sets, introduced by Geelen, Gerards, and Whittle~\cite{GGW09}, can be reduced to a size bounded in terms of their connectivity while preserving their connectivity and titanicity. We also establish conditions under which deletions and contractions within a $\cT$-weak set preserve the order and breadth of the generated tangle. Breadth-criticality and induction on connectivity then yield the required bound for weak sets, from which the main theorem follows.

Section~\ref{sec:preliminaries} gives the definitions and preliminary
results used in the proof. Section~\ref{sec:titanic} first proves a
simultaneous linking result with deletions and contractions restricted
to a specified set, and then uses that result to reduce titanic sets.
Section~\ref{sec:transfer} proves the theorem on generation and
preservation of breadth in minors. Section~\ref{sec:bounds} establishes
the bound for weak sets and completes the proof of
Theorem~\ref{thm:main}.

\section{Preliminaries}
\label{sec:preliminaries}

All sets and matroids considered in this paper are finite. By convention, a partition may have empty parts.
We write $2^E$ for the collection of all subsets of a set $E$, and
$\sqcup$ for disjoint union. A collection of subsets \emph{covers}
$E$ if its union equals $E$. The sets in a cover need not
be disjoint, and repetitions are allowed.

We use the rank axiomatization of a matroid; see
\cite[Chapter~1]{Oxley11}. A \emph{matroid} $M$ on a ground set $E(M)$
is specified by an integer-valued \emph{rank function} $r_M$ on
$2^{E(M)}$ such that, for all $X,Y\subseteq E(M)$,
\begin{enumerate}
\item $0\le r_M(X)\le |X|$;
\item $r_M(X)\le r_M(Y)$ whenever $X\subseteq Y$;
\item $r_M(X)+r_M(Y)\ge r_M(X\cup Y)+r_M(X\cap Y)$.
\end{enumerate}
The rank of $M$ is $r(M)=r_M(E(M))$. A set $X$ is \emph{independent}
if $r_M(X)=|X|$, and \emph{spanning} if $r_M(X)=r(M)$. An independent
spanning set is a \emph{basis}. A set that is not spanning is
\emph{nonspanning}. An element $e$ is a \emph{loop} if
$r_M(\{e\})=0$.

The \emph{closure} of $X\subseteq E(M)$ is
\[
 \cl_M(X)=\{e\in E(M):r_M(X\cup\{e\})=r_M(X)\}.
\]
A set $F$ is a \emph{flat} if $\cl_M(F)=F$; a flat is \emph{proper}
if $F\ne E(M)$. Every proper flat has rank strictly less than $r(M)$, and a flat of rank $r(M)-1$ is a \emph{hyperplane}. Note that
$r_M(\cl_M(X))=r_M(X)$ and
$\cl_M(X)\subseteq\cl_M(Y)$ for $X\subseteq Y$.

For $S\subseteq E(M)$, the \emph{restriction} $M|S$ has ground set
$S$ and rank function $r_{M|S}(X)=r_M(X)$ for all $X\subseteq S$.
The \emph{deletion} $M\setminus D$ is $M|(E(M)\setminus D)$.
For $C\subseteq E(M)$, the \emph{contraction} $M/C$
has ground set $E(M)\setminus C$.
For every $X\subseteq E(M)\setminus C$, its rank
function is given by
\begin{equation}\label{eq:contraction-rank}
r_{M/C}(X)=r_M(X\cup C)-r_M(C).
\end{equation}
A \emph{minor} is a matroid $N=(M/C)\setminus D$, where $C$ and $D$
are disjoint subsets of $E(M)$. A minor is \emph{proper}
if $C\cup D\ne\varnothing$. We write $M/e$ and $M\setminus e$ when
$C=\{e\}$ or $D=\{e\}$, respectively. Deletion and contraction on
distinct elements commute.

The \emph{dual matroid} $M^*$ has ground set $E(M)$ and rank
function given, for every $X\subseteq E(M)$, by
\[
r_{M^*}(X)=|X|+r_M(E(M)\setminus X)-r(M).
\]
The dual rank formula implies that $(M/e)^*=M^*\setminus e$ and
$(M\setminus e)^*=M^*/e$. A matroid $M$ is \emph{representable}
over a field if there is a family of vectors indexed by $E(M)$ for
which a subset of $E(M)$ is independent in $M$ exactly when its
indexed vectors are linearly independent.

The connectivity function $\lambda_M$ defined in
Section~\ref{sec:introduction} is nonnegative and symmetric, that is, 
$\lambda_M(X)=\lambda_M(E(M)\setminus X)$. It satisfies
$\lambda_M(\varnothing)=0$ and $\lambda_M(X)\le r_M(X)\le |X|$. It is also \emph{submodular}, meaning that for all $X,Y\subseteq E(M)$,
\begin{equation}\label{eq:submodular}
 \lambda_M(X)+\lambda_M(Y)
 \ge\lambda_M(X\cup Y)+\lambda_M(X\cap Y).
\end{equation}
Applying~\eqref{eq:submodular} to $X$ and $E(M)\setminus Y$ and
using symmetry yields the \emph{posimodular} inequality
\begin{equation}\label{eq:posimodular}
 \lambda_M(X)+\lambda_M(Y)
 \ge\lambda_M(X\setminus Y)+\lambda_M(Y\setminus X).
\end{equation}
These properties follow from the rank axioms and the dual rank formula gives $\lambda_{M^*}(X)=\lambda_M(X)$.

For a positive integer $k$, a \emph{$k$-separation} of $M$ is a partition
$E(M)=A\sqcup B$ with $|A|,|B|\ge k$ and $\lambda_M(A)\le k-1$.
A matroid is \emph{$k$-connected} if it has no $j$-separation for any
$1\le j<k$. A matroid is \emph{connected} if it is $2$-connected.
A matroid is \emph{weakly $4$-connected} if it is $3$-connected and
$\lambda_M(A)\ge3$ whenever $|A|,|E(M)\setminus A|\ge5$.
Equivalently, a matroid is weakly $4$-connected if and only if it is
$(0,1,4)$-connected in the notation of Section~\ref{sec:introduction}.

For disjoint sets $A,B\subseteq E(M)$, define the
\emph{connectivity between $A$ and $B$} by
\begin{equation}\label{eq:kappa}
 \kappa_M(A,B)=
 \min\{\lambda_M(Z):A\subseteq Z\subseteq E(M)\setminus B\}.
\end{equation}
The minimum is attained because the ground set is finite, and the symmetry of $\lambda_M$ gives $\kappa_M(A,B)=\kappa_M(B,A)$.

The following lemma records two monotonicity
properties of connectivity under taking minors. See~\cite[Lemma~26]{BJOSW25} and~\cite[Lemma~2.4]{HVZ14}.

\begin{lemma}\label{lem:monotonicity}
	Let $N$ be a minor of a matroid $M$.
	For every $X\subseteq E(M)$,
	\begin{equation}\label{eq:trace-monotonicity}
		\lambda_N(X\cap E(N))\le\lambda_M(X).
	\end{equation}
	If $A$ and $B$ are disjoint subsets of $E(N)$, then $0\le\kappa_N(A,B)\le\kappa_M(A,B).$
\end{lemma}
We use the following definition of a matroid tangle from~\cite{BJOSW25}.

\begin{definition}\label{def:tangle}
Let $q\ge2$. A \emph{tangle of order $q$} in a matroid $M$ is a
collection $\cT\subseteq2^{E(M)}$ satisfying the following axioms.
\begin{enumerate}[label=\textup{(T\arabic*)}]
\item $\lambda_M(A)<q-1$ for every $A\in\cT$.
\item If $A\subseteq E(M)$ and $\lambda_M(A)<q-1$, then either $A\in\cT$
  or $E(M)\setminus A\in\cT$.
\item If $A,B,C\in\cT$, then $A\cup B\cup C\ne E(M)$ (where $A,B$, and $C$ need not be distinct).
\item $E(M)\setminus\{e\}\notin\cT$ for every $e\in E(M)$.
\end{enumerate}
The members of $\cT$ are called \emph{$\cT$-small}, and a subset of $E(M)$ is
\emph{$\cT$-weak} if it is contained in a $\cT$-small set.
\end{definition}

\begin{remark}\label{rem:tangle-properties}
Let $\cT$ be a tangle of order $q$ in a matroid $M$.
If $A,B\in\cT$, $X\subseteq A\cup B$, and
$\lambda_M(X)<q-1$, then $X\in\cT$.
Indeed, otherwise $E(M)\setminus X\in\cT$ by~(T2),
and the three sets $A$, $B$, and $E(M)\setminus X$
would cover $E(M)$, contrary to~(T3).
In particular, if $X\subseteq A\in\cT$ and
$\lambda_M(X)<q-1$, then $X\in\cT$.
Moreover, if $q\ge3$, then every singleton belongs to $\cT$: since $\lambda_M(\{e\})\le1<q-1$ for every
$e\in E(M)$, axioms~(T2) and~(T4) imply that
$\{e\}\in\cT$.
\end{remark}

The tangle matroid was introduced in~\cite{GGRW06}.
Let $\cT$ be a tangle of order $q$ in $M$.
By~\cite[Theorem~13 and Lemma~15(i)]{BJOSW25}, the function
\begin{equation}\label{eq:tangle-rank}
\rho_{M,\cT}(X)=
\min\bigl(\{q-1\}\cup
\{\lambda_M(A):X\subseteq A,\ A\in\cT\}\bigr)
\end{equation}
defined for every $X\subseteq E(M)$ is the rank function of a matroid of rank $q-1$
on $E(M)$. This matroid is called the
\emph{tangle matroid} of $\cT$ and is
denoted by $M_{\cT}$. 

A \emph{uniform matroid} of rank $r$ on a set of $n\ge r$ elements,
denoted by $U_{r,n}$, has rank function $r(X)=\min\{r,|X|\}$.
A restriction $L|U$ of a rank-$r$ matroid $L$ is a
\emph{spanning uniform restriction} if $|U|\ge r$ and
$r_L(X)=\min\{r,|X|\}$ for every $X\subseteq U$.

Let $\cT$ be a tangle of order $q$ in $M$. Its \emph{breadth},
denoted by $\br_M(\cT)$, is the maximum size of a spanning uniform
restriction of $M_{\cT}$. A set $U\subseteq E(M)$ is called a
\emph{witness for $\br_M(\cT)$} if $M_{\cT}|U$ is a spanning
uniform restriction of $M_{\cT}$ and $|U|=\br_M(\cT)$.

Let $N$ be a minor of $M$. We say that $\cT$ \emph{generates} a
tangle $\cS$ of order $q$ in $N$ if $\cS$ is the unique tangle of
order $q$ in $N$ satisfying $A\cap E(N)\in\cS$ for every
$A\in\cT$. The tangle $\cT$ is \emph{breadth-critical} in $M$ if,
for every proper minor $N$ in which $\cT$ generates a tangle $\cS$
of order $q$, one has $\br_N(\cS)<\br_M(\cT)$.

The following lemma records consequences of results
of Brettell et al.~\cite{BJOSW25}.
\begin{lemma}\label{lem:critical-connectivity}
Let $\cT$ be a breadth-critical tangle of order
$q\ge2$ in a matroid $M$.
Then $M$ is connected. Moreover, $M$ is
$3$-connected if $q\ge3$, and weakly
$4$-connected if $q\ge4$.
\end{lemma}

\begin{proof}
The first assertion follows from
\cite[Corollary~50]{BJOSW25}.
Suppose that $q\ge3$ and $M$ is not
$3$-connected. By~\cite[Corollary~53]{BJOSW25},
there is a proper minor $N$ of $M$ in which
$\cT$ generates a tangle $\cS$ of order $q$
with $\br_N(\cS)=\br_M(\cT)$.
This contradicts breadth-criticality.
The last assertion follows directly from
\cite[Corollary~66]{BJOSW25}.
\end{proof}

Titanic sets were introduced in~\cite[Section~4]{GGW09}. We adopt the
formulation of~\cite[Section~4]{BJOSW25}.

\begin{definition}\label{def:titanic}
A set $H\subseteq E(M)$ is \emph{titanic} in $M$ if every partition
$H=H_1\sqcup H_2\sqcup H_3$ has at least one part $H_j$ satisfying
$\lambda_M(H_j)\ge\lambda_M(H)$. In particular, every set of
connectivity zero is titanic.
\end{definition}

The following covering formulation is given
in~\cite[Lemma~35]{BJOSW25}.

\begin{lemma}\label{lem:cover}
	Let $H$ be a titanic set in a matroid $M$.
	If $H_1,H_2,H_3\subseteq H$ satisfy
	$H_1\cup H_2\cup H_3=H$, then
	$\lambda_M(H_j)\ge\lambda_M(H)$
	for at least one $j\in\{1,2,3\}$.
\end{lemma}

\section{Reducing titanic sets}
\label{sec:linking}
\label{sec:titanic}

We use the following single-element form of Tutte's linking theorem;
see~\cite{Tutte65} and~\cite[Lemma~2.1]{CW14}.

\begin{theorem}\label{thm:tutte}
Let $M$ be a matroid, let $A$ and $B$ be disjoint subsets of $E(M)$,
and let
$e\in E(M)\setminus(A\cup B)$. Then at least one of $M/e$ and
$M\setminus e$ preserves the connectivity between $A$ and $B$, that is, 
\[
 \kappa_{M/e}(A,B)=\kappa_M(A,B)
 \quad\text{or}\quad
 \kappa_{M\setminus e}(A,B)=\kappa_M(A,B).
\]
\end{theorem}

Chen and Whittle~\cite[Theorem~1.1]{CW14} proved the following
result for two pairs of sets.
\begin{theorem}\label{thm:cw}
Let $M$ be a matroid, and let $A_1,B_1,A_2,B_2\subseteq E(M)$ satisfy
$A_1\cap B_1=A_2\cap B_2=\varnothing$. Put
$k=\kappa_M(A_1,B_1)$ and $\ell=\kappa_M(A_2,B_2)$, and assume
$k,\ell\ge1$. If
\[
 |E(M)\setminus(A_1\cup B_1\cup A_2\cup B_2)|
 \ge(2\ell+1)2^{2k+1},
\]
then there is an element $e\in E(M)\setminus(A_1\cup B_1\cup A_2\cup B_2)$ such that at least one
of the two matroids $M/e$ and $M\setminus e$ preserves both
$\kappa_M(A_1,B_1)$ and $\kappa_M(A_2,B_2)$.
\end{theorem}

The two pairs in Theorem~\ref{thm:cw} may overlap. The following
proposition restricts all deletions and contractions to a specified
set while preserving a finite family of connectivities. 
Recall that $b_t=(2t+1)2^{2t+1}$ from~\eqref{eq:constants}.

\begin{proposition}\label{prop:linking}
Let $M$ be a matroid, and let $m\ge1$ and $t\ge0$ be integers.
For every $j\in\{1,\ldots,m\}$, let
$A_j,B_j\subseteq E(M)$ be disjoint and satisfy
$\kappa_M(A_j,B_j)\le t$. Let $P\subseteq E(M)\setminus\bigcup_{j=1}^{m}(A_j\cup B_j)$.
Then there exist disjoint subsets $C,D\subseteq P$ such that the minor
$N=(M/C)\setminus D$ satisfies
\[
 |P\cap E(N)|\le\binom m2(b_t-1)
\]
and $\kappa_N(A_j,B_j)=\kappa_M(A_j,B_j)$ for every
$j\in\{1,\ldots,m\}$. 
\end{proposition}

\begin{proof}
Write $k_j=\kappa_M(A_j,B_j)$ for each $j\in\{1,\ldots,m\}$. Choose disjoint $C,D\subseteq P$
so that $N=(M/C)\setminus D$ preserves all the values $k_j$ and
$|C\cup D|$ is as large as possible. Such a choice exists because
$C=D=\varnothing$ is admissible and $P$ is finite. Put
$F=P\setminus(C\cup D)=P\cap E(N)$.

For each $e\in F$ and $j\in\{1,\ldots,m\}$, let $O_j(e)$ be the set of operations in
$\{/,\setminus\}$ that preserve $k_j$ in $N$, that is,
\[
O_j(e)=\{\circ\in\{/,\setminus\}:
\kappa_{N\circ e}(A_j,B_j)=k_j\},
\]
where $N\circ e$ denotes $N/e$ when $\circ=/$, and
$N\setminus e$ when $\circ=\setminus$.
Theorem~\ref{thm:tutte} implies that every
$O_j(e)$ is nonempty. For every $e\in F$, the maximality of
$|C\cup D|$ implies that
\begin{equation}\label{eq:no-common-operation}
\bigcap_{j=1}^{m}O_j(e)=\varnothing.
\end{equation}

If $k_j=0$, Lemma~\ref{lem:monotonicity} shows that
$O_j(e)=\{/,\setminus\}$ for every $e\in F$.
If $t=0$, then $k_j=0$ for every $j$, so
\eqref{eq:no-common-operation} implies $F=\varnothing$.
If $m=1$, equation~\eqref{eq:no-common-operation} contradicts
nonemptiness of $O_1(e)$ for every surviving $e$, so
$F=\varnothing$. We may therefore assume $m\ge2$ and $t\ge1$.

For $1\le i<j\le m$, define
\[
F_{ij}=\{e\in F:O_i(e)\cap O_j(e)=\varnothing\}.
\]
Fix $1\le i<j\le m$. We first show that
$|F_{ij}|\le b_t-1$.
If $k_i=0$ or $k_j=0$, then $F_{ij}=\varnothing$.
We may therefore assume that $k_i,k_j\ge1$.

Suppose, for a contradiction, that $|F_{ij}|\ge b_t$.
Set $N_0=N$. We apply Theorem~\ref{thm:cw}
to the two pairs $(A_i,B_i)$ and
$(A_j,B_j)$ in $N_0$.
The two connectivities in $N_0$ are $k_i$ and $k_j$,
respectively, and $F_{ij}$ is disjoint from
$A_i\cup B_i\cup A_j\cup B_j$. Since
\[
|F_{ij}|\ge b_t=(2t+1)2^{2t+1}
\ge(2k_j+1)2^{2k_i+1},
\]
Theorem~\ref{thm:cw} gives an element
$e_0\in E(N_0)\setminus
(A_i\cup B_i\cup A_j\cup B_j)$ and
an operation $\circ_0\in\{/,\setminus\}$
such that $\kappa_{N_0\circ_0 e_0}(A_v,B_v)=k_v$ for $v=i,j$.

If $e_0\notin F_{ij}$, put
$N_1=N_0\circ_0 e_0$ and apply
Theorem~\ref{thm:cw} again to
$(A_i,B_i)$ and $(A_j,B_j)$ in $N_1$.
More generally, whenever the theorem gives
an element $e_h\notin F_{ij}$, put
$N_{h+1}=N_h\circ_h e_h$ and apply
Theorem~\ref{thm:cw} to the two pairs
$(A_i,B_i)$ and $(A_j,B_j)$ in $N_{h+1}$.

Each such operation preserves both connectivities
and retains every element of $F_{ij}$. Since
$|F_{ij}|\ge b_t$, Theorem~\ref{thm:cw}
remains applicable after each such operation.
As $N$ is finite, some application must supply
an element $e_h\in F_{ij}$.

Let $K=N_h$, $e=e_h$, and $\circ=\circ_h$.
For each $v\in\{i,j\}$, the choice of $e$ and $\circ$ gives
$\kappa_{K\circ e}(A_v,B_v)=k_v$.
All elements deleted or contracted in obtaining
$K$ from $N$ lie outside
$F_{ij}\cup A_i\cup B_i\cup A_j\cup B_j$.
Thus we may write
\[
K=N/C'\setminus D',
\]
where $C'\cap D'=\varnothing$,
$e\notin C'\cup D'$, and
$C'\cup D'$ is disjoint from
$A_i\cup B_i\cup A_j\cup B_j$.

Deletion and contraction on distinct
elements commute. Hence
\[
K\circ e=(N\circ e)/C'\setminus D'.
\]
In particular, $K\circ e$ is a minor of
$N\circ e$, and
$A_i\cup B_i\cup A_j\cup B_j
\subseteq E(K\circ e)\subseteq E(N\circ e)$.
Thus, for $v\in\{i,j\}$,
Lemma~\ref{lem:monotonicity} gives
\[
\begin{aligned}
k_v
&=\kappa_N(A_v,B_v)\\
&\ge\kappa_{N\circ e}(A_v,B_v)\\
&\ge\kappa_{K\circ e}(A_v,B_v)=k_v.
\end{aligned}
\]
Consequently, $\kappa_{N\circ e}(A_v,B_v)=k_v$
for each $v\in\{i,j\}$.
By the definition of $O_i(e)$ and
$O_j(e)$, we have
$\circ\in O_i(e)\cap O_j(e)$.
This contradicts $e\in F_{ij}$.
Therefore $|F_{ij}|\le b_t-1$.

Finally, let $e\in F$.
Since each $O_j(e)$ is a nonempty subset of
$\{/,\setminus\}$,
equation~\eqref{eq:no-common-operation}
implies that there are indices $i<j$ such that
$O_i(e)\cap O_j(e)=\varnothing$.
By the definition of $F_{ij}$, we have
$e\in F_{ij}$.
Since $F_{ij}\subseteq F$ for every $i<j$,
it follows that
\[
F=\bigcup_{1\le i<j\le m}F_{ij}.
\]
Therefore,
\[
|F|
\le\sum_{1\le i<j\le m}|F_{ij}|
\le\binom{m}{2}(b_t-1).
\]
\end{proof}

The following result is due to Geelen, Gerards,
and Whittle~\cite[Lemma~4.6]{GGW07}.

\begin{lemma}\label{lem:linkage-matroid}
Let $M$ be a matroid, let $H\subseteq E(M)$,
and put $H^c=E(M)\setminus H$.
There exists a matroid $L$ with ground set $H$
whose rank function satisfies
$r_L(X)=\kappa_M(X,H^c)$ for every $X\subseteq H$.
Moreover, $r(L)=\lambda_M(H)$.
\end{lemma}

\begin{lemma}\label{lem:certificate}
Let $L$ be a matroid of rank $t$ whose ground set
cannot be covered by three proper flats, not necessarily distinct.
Then there exists a spanning set $W\subseteq E(L)$ with
$|W|\le c_t=(3^t-1)/2$ such that no three proper flats of
$L|W$ cover $W$.
\end{lemma}

\begin{proof}
If $t=0$, the only flat of $L$ is $E(L)$, so $L$ has no proper
flats. The restriction to $W=\varnothing$ also has rank zero
and no proper flats. This proves the rank-zero case.

Suppose that $t\ge1$. Since no three proper flats of $L$ cover
$E(L)$, choose an inclusion-minimal set $W\subseteq E(L)$ such that
$W\not\subseteq J_1\cup J_2\cup J_3$ for every three proper flats
$J_1,J_2,J_3$ of $L$.
If $W$ does not span $L$, then $\cl_L(W)$ is a proper flat of $L$.
Taking $J_1=J_2=J_3=\cl_L(W)$ gives
$W\subseteq J_1\cup J_2\cup J_3$, contrary to the choice of $W$.
Thus $W$ spans $L$.

Suppose that three proper flats $J_1,J_2,J_3$
of $L|W$ cover $W$. Since $W$ spans $L$,
none of $J_1,J_2,J_3$ spans $L$.
Thus $\cl_L(J_1),\cl_L(J_2),\cl_L(J_3)$
are proper flats of $L$, and
\[
W\subseteq
\cl_L(J_1)\cup\cl_L(J_2)\cup\cl_L(J_3),
\]
contrary to the choice of $W$.

For each $e\in W$, the minimality of $W$
yields three proper flats
$F_{e,1},F_{e,2},F_{e,3}$ of $L$ such that
\[
W\setminus\{e\}
\subseteq F_{e,1}\cup F_{e,2}\cup F_{e,3}.
\]
For each $j\in\{1,2,3\}$, we have $e\notin F_{e,j}$;
otherwise $F_{e,1}\cup F_{e,2}\cup F_{e,3}$ would contain $W$.
Consequently, $W$ contains no loops of $L$,
since every flat contains every loop.

We prove by induction on $d$ that every flat $J$
of $L$ with $r_L(J)=d$ satisfies
\[
|W\cap J|\le c_d.
\]
If $d=0$, then $J$ consists of loops, so
$W\cap J=\varnothing$ and $|W\cap J|=0=c_0$.

Suppose that $d\ge1$ and that the inequality
holds for all flats of rank less than $d$.
Let $J$ be a flat of $L$ of rank $d$.
If $W\cap J=\varnothing$, the inequality holds.
Otherwise, choose $e\in W\cap J$.
For each $j\in\{1,2,3\}$,
the set $J\cap F_{e,j}$ is a proper flat
of $L$ contained in $J$, since
$e\in J$ but $e\notin F_{e,j}$.
Hence $J\cap F_{e,j}$ has rank less than $d$.
By the induction hypothesis and the monotonicity of $(c_d)$,
\[
|W\cap J\cap F_{e,j}|\le c_{d-1}.
\]
Since
\[
(W\cap J)\setminus\{e\}
\subseteq
\bigcup_{j=1}^{3}(W\cap J\cap F_{e,j}),
\]
we obtain
\[
\begin{aligned}
|W\cap J|
&\le1+\sum_{j=1}^{3}|W\cap J\cap F_{e,j}|\\
&\le1+3c_{d-1}=c_d.
\end{aligned}
\]
This completes the induction.
Taking $J=E(L)$ gives $|W|\le c_t$.
\end{proof}

Recall from \eqref{eq:constants} that $h_t=c_t+\binom{2^{c_t}}2(b_t-1)$.

\begin{proposition}\label{prop:compression}
Let $M$ be a matroid, let $H\subseteq E(M)$ be
titanic, and put $t=\lambda_M(H)$.
There exist disjoint sets $C,D\subseteq H$ such that,
for $N=(M/C)\setminus D$ and
$H'=H\setminus(C\cup D)$,
\[
|H'|\le h_t,\qquad \lambda_N(H')=t,
\]
and $H'$ is titanic in $N$.
\end{proposition}

\begin{proof}
If $t=0$, take $C=\varnothing$ and $D=H$.
Then $H'=\varnothing$, so
$|H'|=h_0=0$ and $\lambda_N(H')=0$.
Moreover, $\varnothing$ is titanic in $N$.

Suppose that $t\ge1$, and put $H^c=E(M)\setminus H$.
By Lemma~\ref{lem:linkage-matroid}, there is a matroid
$L$ with ground set $H$ whose rank function satisfies
$r_L(X)=\kappa_M(X,H^c)$ for every $X\subseteq H$.
Moreover, $r(L)=\lambda_M(H)=t$.

We first show that no three proper flats of $L$
cover $H$.
Suppose that $F_1,F_2,F_3$ are proper flats of
$L$ whose union is $H$.
For each $j\in\{1,2,3\}$, the definition of $r_L$ gives
\[
r_L(F_j)=\kappa_M(F_j,H^c)
=\min_{F_j\subseteq Z\subseteq H}\lambda_M(Z).
\]
Since $F_j$ is a proper flat of $L$, we have $r_L(F_j)<t$.

For each $j\in\{1,2,3\}$, choose
$Z_j$ with $F_j\subseteq Z_j\subseteq H$
such that $\lambda_M(Z_j)=r_L(F_j)<t$.
Since $F_1\cup F_2\cup F_3=H$, we also have
$Z_1\cup Z_2\cup Z_3=H$.
This contradicts Lemma~\ref{lem:cover},
because $H$ is titanic and
$\lambda_M(Z_j)<\lambda_M(H)$
for every $j\in\{1,2,3\}$.

By Lemma~\ref{lem:certificate}, there is a set
$W\subseteq H$ that spans $L$, satisfies
$|W|\le c_t$, and cannot be covered by
three proper flats of $L|W$.

Put $P=H\setminus W$. For every $X\subseteq W$,
consider the pair $(X,H^c)$. There are $2^{|W|}$
such pairs, and
$\kappa_M(X,H^c)=r_L(X)\le t$.
Since $P\cap(W\cup H^c)=\varnothing$,
Proposition~\ref{prop:linking} gives disjoint
sets $C,D\subseteq P$ such that
$N=(M/C)\setminus D$ satisfies
\begin{equation}\label{eq:certificate-preserved}
\kappa_N(X,H^c)=r_L(X)
\end{equation}
for every $X\subseteq W$,
and
\[
|P\setminus(C\cup D)|
\le\binom{2^{|W|}}{2}(b_t-1).
\]

Put $H'=H\setminus(C\cup D)$.
Since $H=W\sqcup P$ and $C,D\subseteq P$,
we have $H'=W\sqcup(P\setminus(C\cup D))$.
Therefore,
\[
\begin{aligned}
|H'|
&=|W|+|P\setminus(C\cup D)|\\
&\le |W|+\binom{2^{|W|}}{2}(b_t-1)\\
&\le c_t+\binom{2^{c_t}}{2}(b_t-1)
=h_t.
\end{aligned}
\]

Since $W$ spans $L$, we have $r_L(W)=t$.
Equation~\eqref{eq:certificate-preserved}
then gives $\kappa_N(W,H^c)=t$.
Moreover, $W\subseteq H'=E(N)\setminus H^c$.
By the definition of $\kappa_N$ and
Lemma~\ref{lem:monotonicity},
\[
t=\kappa_N(W,H^c)
\le\lambda_N(H')
\le\lambda_M(H)=t.
\]
Hence $\lambda_N(H')=t$.

It remains to show that $H'$ is titanic in $N$.
Consider a partition $H'=J_1\sqcup J_2\sqcup J_3$.
For each $j\in\{1,2,3\}$, put $X_j=W\cap J_j$.
Then $W=X_1\cup X_2\cup X_3$.

At least one of $X_1,X_2,X_3$ spans $L$.
Otherwise, the three sets
$\cl_{L|W}(X_1)$, $\cl_{L|W}(X_2)$,
and $\cl_{L|W}(X_3)$ would be proper flats
of $L|W$ covering $W$, contrary to the
choice of $W$.

Choose $j\in\{1,2,3\}$ such that $X_j$ spans $L$.
Then $r_L(X_j)=t$.
Since $X_j\subseteq J_j\subseteq E(N)\setminus H^c$,
equation~\eqref{eq:certificate-preserved} gives
\[
\lambda_N(J_j)
\ge\kappa_N(X_j,H^c)
=r_L(X_j)=t=\lambda_N(H').
\]
Thus $H'$ is titanic in $N$.
\end{proof}

\section{Tangles in minors}
\label{sec:transfer}
\begin{lemma}\label{lem:boundary}
Let $M$ be a matroid, let $H\subseteq E(M)$, and put
$H^c=E(M)\setminus H$. Let $C,D\subseteq H$ be disjoint,
and set $N=(M/C)\setminus D$ and $H'=H\setminus(C\cup D)$.
If $\lambda_N(H')=\lambda_M(H)$, then, for every $Z\subseteq H^c$,
\begin{equation}\label{eq:exterior}
\lambda_N(Z)=\lambda_M(Z)
\qquad \text{and} \qquad
\lambda_N(Z\cup H')=\lambda_M(Z\cup H).
\end{equation}
\end{lemma}

\begin{proof}
If $C\cup D=\varnothing$, then $N=M$ and $H'=H$,
so both equalities hold. Suppose that $C\cup D\ne\varnothing$.
Put $n=|C\cup D|$, and order the elements of $C\cup D$ as
$e_1,\ldots,e_n$. Starting with $M_0=M$, define $M_a=M_{a-1}/e_a$
if $e_a\in C$, and $M_a=M_{a-1}\setminus e_a$ if $e_a\in D$,
for each $a\in\{1,\ldots,n\}$. Then $M_n=N$.
For each $a\in\{0,\ldots,n\}$, put $H_a=H\cap E(M_a)$.
Since $C,D\subseteq H$, we have $E(M_a)=H_a\sqcup H^c$.
Lemma~\ref{lem:monotonicity} gives
\[
\lambda_M(H)=\lambda_{M_0}(H_0)
\ge\lambda_{M_1}(H_1)\ge\cdots
\ge\lambda_{M_n}(H_n)=\lambda_N(H')=\lambda_M(H).
\]
Thus every inequality is an equality. By symmetry of connectivity,
$\lambda_{M_a}(H^c)=\lambda_{M_{a-1}}(H^c)$ for every
$a\in\{1,\ldots,n\}$.

Consider a step that contracts $e_a$, and put $K=M_{a-1}$,
$e=e_a$, and $\varepsilon=r_K(\{e\})$.
Since $e\notin H^c$,
the definition of $\lambda$ and equation~\eqref{eq:contraction-rank} give
\begin{align*}
\lambda_{K/e}(H^c)
&=r_{K/e}(H^c)
  +r_{K/e}(E(K/e)\setminus H^c)-r(K/e)\\
&=r_K(H^c\cup\{e\})+r_K(E(K)\setminus H^c)
  -r(K)-\varepsilon,
\end{align*}
and
\[\lambda_K(H^c)=r_K(H^c)+r_K(E(K)\setminus H^c)-r(K).\]
Since $\lambda_{K/e}(H^c)=\lambda_K(H^c)$, we have 
\[
r_K(H^c\cup\{e\})-r_K(H^c)=\varepsilon.
\]
Let $Z\subseteq H^c$. Submodularity of $r_K$, applied to
$Z\cup\{e\}$ and $H^c$, gives
\[
r_K(Z\cup\{e\})-r_K(Z)
\ge r_K(H^c\cup\{e\})-r_K(H^c)=\varepsilon.
\]
Applying submodularity to $Z$ and $\{e\}$ gives the reverse
inequality. Hence $r_K(Z\cup\{e\})-r_K(Z)=\varepsilon$.
Since $e\notin Z$, the definition of $\lambda$ and
equation~\eqref{eq:contraction-rank} give
\begin{align*}
\lambda_{K/e}(Z)
&=r_{K/e}(Z)
  +r_{K/e}(E(K/e)\setminus Z)-r(K/e)\\
&=r_K(Z\cup\{e\})+r_K(E(K)\setminus Z)
  -r(K)-\varepsilon\\
&=r_K(Z)+r_K(E(K)\setminus Z)-r(K)\\
&=\lambda_K(Z).
\end{align*}
The penultimate equality uses
$r_K(Z\cup\{e\})-r_K(Z)=\varepsilon$.

For a step that deletes $e$ from $K$, use
$(K\setminus e)^*=K^*/e$ and $\lambda_{K^*}=\lambda_K$.
The equality $\lambda_{K\setminus e}(H^c)=\lambda_K(H^c)$
becomes $\lambda_{K^*/e}(H^c)=\lambda_{K^*}(H^c)$.
The contraction argument applied to $K^*$ therefore gives
$\lambda_{K\setminus e}(Z)=\lambda_K(Z)$ for every $Z\subseteq H^c$.

It follows that $\lambda_{M_a}(Z)=\lambda_{M_{a-1}}(Z)$ for
every $a\in\{1,\ldots,n\}$ and every $Z\subseteq H^c$.
Thus $\lambda_N(Z)=\lambda_M(Z)$ for every $Z\subseteq H^c$.
Finally, $E(N)=H'\sqcup H^c$ and $E(M)=H\sqcup H^c$.
For every $Z\subseteq H^c$, symmetry gives
\[
\lambda_N(Z\cup H')
=\lambda_N(H^c\setminus Z)
=\lambda_M(H^c\setminus Z)
=\lambda_M(Z\cup H).
\]
\end{proof}

The following theorem gives sufficient conditions for a tangle to generate
a tangle of the same order in a minor and for their breadths to be equal.

\begin{theorem}\label{thm:transfer}
Let $\cT$ be a tangle of order $q\ge3$ in a matroid $M$,
and let $H\subseteq E(M)$ be $\cT$-weak.
Let $C,D\subseteq H$ be disjoint, and put
$N=(M/C)\setminus D$ and $H'=H\setminus(C\cup D)$.
Suppose that $H'$ is titanic in $N$ and
$\lambda_N(H')=\lambda_M(H)$.
Then $\cT$ generates a tangle $\cS$ of order $q$ in $N$.
If there is a witness $U$ for $\br_M(\cT)$ with
$U\cap H=\varnothing$, then $\br_N(\cS)=\br_M(\cT)$.
\end{theorem}

\begin{proof}
Put $k=q-1$, $t=\lambda_M(H)$, and
$H^c=E(M)\setminus H$.
Since $H$ is $\cT$-weak, choose $F\in\cT$ with $H\subseteq F$,
and put $F'=F\cap E(N)$.
We have $E(N)=H'\sqcup H^c$ and
$E(M)\setminus E(N)=C\cup D\subseteq H\subseteq F$.
Put
\[
\cW=\{B\subseteq E(N):
B\subseteq A\cap E(N)\text{ for some }A\in\cT\}.
\]

For every partition $E(N)=B_1\sqcup B_2$,
at most one of $B_1,B_2$ belongs to $\cW$.
Otherwise, there would be $A_1,A_2\in\cT$ with
$B_1\subseteq A_1\cap E(N)$ and
$B_2\subseteq A_2\cap E(N)$.
Since $E(N)\subseteq A_1\cup A_2$ and
$E(M)\setminus E(N)\subseteq F$, we would have
$A_1\cup A_2\cup F=E(M)$, contrary to~(T3).

We construct $\cS$ by selecting one part of
each partition $E(N)=B_1\sqcup B_2$ with
$\lambda_N(B_1)<k$.
If $B_1\in\cW$, select $B_1$; if $B_2\in\cW$,
select $B_2$. As shown above, at most one
part belongs to $\cW$.

Suppose that neither $B_1$ nor $B_2$
belongs to $\cW$.
Since $H'$ is titanic in $N$ and $\lambda_N(H')=t$,
Lemma~\ref{lem:cover}, applied to
$H'=(B_1\cap H')\sqcup(B_2\cap H')$, gives
$\lambda_N(B_1\cap H')\ge t$ or
$\lambda_N(B_2\cap H')\ge t$.
Interchanging $B_1$ and $B_2$ if necessary, assume that
$\lambda_N(B_1\cap H')\ge t$.
Since $B_1\sqcup B_2=E(N)=H'\sqcup H^c$, we have
\[
E(N)\setminus(B_2\cap H^c)=B_1\cup H'.
\]
By Lemma~\ref{lem:boundary} and symmetry of
the connectivity function,
\[
\lambda_M(B_2\cap H^c)
=\lambda_N(B_2\cap H^c)
=\lambda_N(B_1\cup H').
\]
Submodularity applied to $B_1$ and $H'$ gives
\[
\lambda_N(B_1\cup H')
\le\lambda_N(B_1)+\lambda_N(H')
-\lambda_N(B_1\cap H').
\]
Since $\lambda_N(H')=t$ and
$\lambda_N(B_1\cap H')\ge t$, we obtain
\[
\lambda_M(B_2\cap H^c)
\le\lambda_N(B_1)<k.
\]
Since $E(M)=H\sqcup H^c$ and
$E(N)=B_1\sqcup B_2$ with $H^c\subseteq E(N)$,
we have $E(M)\setminus(B_2\cap H^c)=H\cup(B_1\cap H^c)$.
If $H\cup(B_1\cap H^c)$ belonged to $\cT$, then
$B_1\subseteq (H\cup(B_1\cap H^c))\cap E(N)$,
contrary to $B_1\notin\cW$.
Thus $H\cup(B_1\cap H^c)\notin\cT$.
Since $\lambda_M(B_2\cap H^c)<k$, axiom~(T2) gives
$B_2\cap H^c\in\cT$.
Suppose also that $\lambda_N(B_2\cap H')\ge t$.
Interchanging $B_1$ and $B_2$ in the preceding argument
gives $B_1\cap H^c\in\cT$.
Then $(B_1\cap H^c)\cup(B_2\cap H^c)\cup F=E(M)$,
contrary to~(T3).

Hence $\lambda_N(B_2\cap H')<t\le\lambda_N(B_1\cap H')$,
and we select $B_2$.
Let $\cS$ be the collection of the sets selected in this way.

For every $B\in\cS\setminus\cW$, the choice
of $B$ in the preceding construction gives
\begin{equation}\label{eq:new-side}
B\cap H^c\in\cT,\qquad
\lambda_N(B\cap H')<t.
\end{equation}

By construction, every member of $\cS$ has connectivity
less than $k$, and at least one part of every partition
of $E(N)$ with connectivity less than $k$ belongs to
$\cS$. Thus $\cS$ satisfies~(T1) and~(T2).

To verify~(T3), suppose that $B_1,B_2,B_3\in\cS$
and $B_1\cup B_2\cup B_3=E(N)$.
For each $j\in\{1,2,3\}$ with $B_j\in\cW$, choose
$A_j\in\cT$ such that $B_j\subseteq A_j\cap E(N)$.

For each $j\in\{1,2,3\}$ with $B_j\notin\cW$, put
$A_j=B_j\cap H^c$. Then $A_j\in\cT$ by~\eqref{eq:new-side}.
Thus $B_j\cap H^c\subseteq A_j$ for every $j\in\{1,2,3\}$,
and hence $H^c\subseteq A_1\cup A_2\cup A_3$.

Suppose that $B_j\in\cW$ and
$\lambda_M(A_j\cap H)\ge t$ for some $j\in\{1,2,3\}$.
Submodularity gives
\[
\lambda_M(A_j\cup H)
\le\lambda_M(A_j)+\lambda_M(H)-\lambda_M(A_j\cap H)
\le\lambda_M(A_j)<k.
\]
Since $A_j,F\in\cT$ and $A_j\cup H\subseteq A_j\cup F$,
Remark~\ref{rem:tangle-properties} gives $A_j\cup H\in\cT$.
The set $A_j\cup H$, together with the other two sets
among $A_1,A_2,A_3$, covers $E(M)$, contrary to~(T3).
Thus, for each
$j\in\{1,2,3\}$ with $B_j\in\cW$, $\lambda_M(A_j\cap H)<t$.
Lemma~\ref{lem:monotonicity} then gives
$\lambda_N(A_j\cap H')\le\lambda_M(A_j\cap H)<t$.

For each $j\in\{1,2,3\}$, put
$Y_j=A_j\cap H'$ if $B_j\in\cW$, and
$Y_j=B_j\cap H'$ otherwise.
If $B_j\in\cW$, then
$B_j\subseteq A_j\cap E(N)$, so
$B_j\cap H'\subseteq Y_j\subseteq H'$.
The preceding inequality gives
$\lambda_N(Y_j)<t$.
If $B_j\notin\cW$, then
$Y_j=B_j\cap H'$, and
$\lambda_N(Y_j)<t$ follows
from~\eqref{eq:new-side}.
Thus $B_j\cap H'\subseteq Y_j\subseteq H'$
and $\lambda_N(Y_j)<t$ for every
$j\in\{1,2,3\}$.
Since $B_1\cup B_2\cup B_3=E(N)$,
we have $Y_1\cup Y_2\cup Y_3=H'$.
This contradicts Lemma~\ref{lem:cover}, because $H'$ is
titanic in $N$ and $\lambda_N(H')=t$.
Therefore~(T3) holds.

For~(T4), let $e\in E(N)$.
Since $q\ge3$, Remark~\ref{rem:tangle-properties} gives
$\{e\}\in\cT$.
Thus $\{e\}\in\cW$ and $\lambda_N(\{e\})\le1<k$.
The selection rule puts $\{e\}$ in $\cS$ and excludes
$E(N)\setminus\{e\}$.
Hence~(T4) holds, and $\cS$ is a tangle of order $q$ in $N$.

We next prove that $\cT$ generates $\cS$.
For each $A\in\cT$, Lemma~\ref{lem:monotonicity} gives
$\lambda_N(A\cap E(N))\le\lambda_M(A)<k$.
Since $A\cap E(N)\in\cW$, the selection rule gives
$A\cap E(N)\in\cS$. In particular, $F'=F\cap E(N)\in\cS$.

Let $\cS_0$ be a tangle of order $q$ in $N$ such that
$A\cap E(N)\in\cS_0$ for every $A\in\cT$.
If $B\in\cS\cap\cW$, choose $A\in\cT$ with
$B\subseteq A\cap E(N)$.
Since $A\cap E(N)\in\cS_0$ and $\lambda_N(B)<k$,
Remark~\ref{rem:tangle-properties} gives $B\in\cS_0$.
If $B\in\cS\setminus\cW$, equation~\eqref{eq:new-side}
gives $B\cap H^c\in\cT$.
Since $B\cap H^c\subseteq E(N)$, we have $B\cap H^c\in\cS_0$.
Since $H\subseteq F$, we have $H'\subseteq F'$.
As $E(N)=H'\sqcup H^c$, it follows that
$B\subseteq(B\cap H^c)\cup F'$.
Both $B\cap H^c$ and $F'=F\cap E(N)$ belong to
$\cS_0$, and $\lambda_N(B)<k$.
By Remark~\ref{rem:tangle-properties},
$B\in\cS_0$.
Hence $\cS\subseteq\cS_0$.
If $B\in\cS_0\setminus\cS$, then $\lambda_N(B)<k$,
and the selection rule for $\cS$ gives
$E(N)\setminus B\in\cS\subseteq\cS_0$.
This contradicts~(T3) for $\cS_0$.
Therefore $\cS_0=\cS$ and $\cT$ generates $\cS$.

We will also use the following consequence of generation:
for every $A\subseteq E(M)$ with $\lambda_M(A)<k$,
\begin{equation}\label{eq:lift}
A\in\cT\quad\Longleftrightarrow\quad A\cap E(N)\in\cS.
\end{equation}
The forward implication was proved above.
If $A\notin\cT$ and $\lambda_M(A)<k$, axiom~(T2)
gives $E(M)\setminus A\in\cT$.
Thus $E(N)\setminus(A\cap E(N))\in\cS$.
Since $\cS$ cannot contain both a set and its complement,
$A\cap E(N)\notin\cS$. This proves the reverse implication.

To compare the breadths, recall that
$\rho_{M,\cT}$ and $\rho_{N,\cS}$ are the rank
functions of $M_{\cT}$ and $N_{\cS}$, respectively.
We first show that, for every $X\subseteq H^c$,
\begin{equation}\label{eq:rank-outside}
\rho_{N,\cS}(X)=\rho_{M,\cT}(X).
\end{equation}

Let $X\subseteq H^c$. We first prove that
$\rho_{M,\cT}(X)\le\rho_{N,\cS}(X)$.
If $\rho_{N,\cS}(X)=k$, this follows from
$\rho_{M,\cT}(X)\le k$.
Otherwise, by~\eqref{eq:tangle-rank}, there is
a set $Y\in\cS$ containing $X$ such that
$\lambda_N(Y)=\rho_{N,\cS}(X)<k$.

Suppose that $\lambda_N(Y\cap H')\ge t$.
By submodularity and $\lambda_N(H')=t$,
\[
\lambda_N(Y\cup H')
\le\lambda_N(Y)+\lambda_N(H')
-\lambda_N(Y\cap H')
\le \lambda_N(Y)<k.
\]
Since $Y,F'\in\cS$ and
$Y\cup H'\subseteq Y\cup F'$,
Remark~\ref{rem:tangle-properties} gives
$Y\cup H'\in\cS$.
Put $A=H\cup(Y\cap H^c)$.
Then $A\cap E(N)=H'\cup(Y\cap H^c)=Y\cup H'$.
By Lemma~\ref{lem:boundary},
$\lambda_M(A)=\lambda_N(A\cap E(N))=\lambda_N(Y\cup H') \le\lambda_N(Y)<k$.
Equation~\eqref{eq:lift} gives $A\in\cT$.
Since $X\subseteq A$, we obtain
$\rho_{M,\cT}(X)\le \lambda_M(A)\le \lambda_N(Y)=\rho_{N,\cS}(X)$.

Now suppose that $\lambda_N(Y\cap H')<t$.
Since $H'$ is titanic, Lemma~\ref{lem:cover}
applied to the partition
$H'=(Y\cap H')\sqcup(H'\setminus Y)$ gives
$\lambda_N(H'\setminus Y)\ge t$.
Since $Y\setminus H'=Y\cap H^c$, applying
\eqref{eq:posimodular} to $Y$ and $H'$ gives
\[
\lambda_N(Y\cap H^c)
\le\lambda_N(Y)+\lambda_N(H')
-\lambda_N(H'\setminus Y)
\le \lambda_N(Y)<k.
\]
Since $Y\cap H^c\subseteq Y\in\cS$,
Remark~\ref{rem:tangle-properties} gives
$Y\cap H^c\in\cS$.
By Lemma~\ref{lem:boundary},
$\lambda_M(Y\cap H^c)=\lambda_N(Y\cap H^c)
\le \lambda_N(Y)<k$.
Equation~\eqref{eq:lift} gives
$Y\cap H^c\in\cT$.
Since $X\subseteq Y\cap H^c$, we again obtain
$\rho_{M,\cT}(X)\le \lambda_N(Y)=\rho_{N,\cS}(X)$.

We have therefore proved
$\rho_{M,\cT}(X)\le\rho_{N,\cS}(X)$ for every $X\subseteq H^c$.

For the reverse inequality, let $X\subseteq E(N)$.
If $\rho_{M,\cT}(X)=k$, then
$\rho_{N,\cS}(X)\le k=\rho_{M,\cT}(X)$.
Otherwise, by~\eqref{eq:tangle-rank}, there is
a set $A\in\cT$ containing $X$ such that
$\lambda_M(A)=\rho_{M,\cT}(X)$.
Since $A\cap E(N)\in\cS$ and contains $X$,
\[
\rho_{N,\cS}(X)
\le\lambda_N(A\cap E(N))
\le\lambda_M(A)=\rho_{M,\cT}(X).
\]
Hence
\begin{equation}\label{eq:rank-domination}
\rho_{N,\cS}(X)\le\rho_{M,\cT}(X)
\end{equation}
for every $X\subseteq E(N)$.
In particular, the reverse inequality holds
for every $X\subseteq H^c$, proving
\eqref{eq:rank-outside}.

Finally, suppose that $U$ is a witness for
$\br_M(\cT)$ with $U\cap H=\varnothing$.
Recall that both $M_{\cT}$ and $N_{\cS}$
have rank $k$. By the definition of breadth,
$|U|=\br_M(\cT)\ge k$, and
$\rho_{M,\cT}(X)=\min\{k,|X|\}$
for every $X\subseteq U$.

Since $U\subseteq H^c\subseteq E(N)$,
equation~\eqref{eq:rank-outside} gives
\[
\rho_{N,\cS}(X)
=\rho_{M,\cT}(X)
=\min\{k,|X|\}
\]
for every $X\subseteq U$.
Thus $N_{\cS}|U$ is a spanning uniform
restriction of $N_{\cS}$, and
$\br_N(\cS)\ge|U|=\br_M(\cT)$.

For the reverse inequality, let $V$ be a
witness for $\br_N(\cS)$. By definition,
$|V|=\br_N(\cS)\ge k$, and
$\rho_{N,\cS}(X)=\min\{k,|X|\}$
for every $X\subseteq V$.

Since $M_{\cT}$ has rank $k$, its rank
function satisfies
$\rho_{M,\cT}(X)\le\min\{k,|X|\}$.
Together with~\eqref{eq:rank-domination},
this gives
\[
\min\{k,|X|\}
=\rho_{N,\cS}(X)
\le\rho_{M,\cT}(X)
\le\min\{k,|X|\}
\]
for every $X\subseteq V$.
Hence $\rho_{M,\cT}(X)=\min\{k,|X|\}$
for every $X\subseteq V$.
As $|V|\ge k$, the restriction
$M_{\cT}|V$ is a spanning uniform restriction of $M_{\cT}$.
Therefore
$|V|\le\br_M(\cT)$, and hence
$\br_N(\cS)\le\br_M(\cT)$.
Combining the two inequalities, we obtain
$\br_N(\cS)=\br_M(\cT)$.
\end{proof}

\section{Separation bounds}
\label{sec:bounds}

We first bound the size of a $\cT$-weak set disjoint from a
witness for $\br_M(\cT)$.
\begin{proposition}\label{prop:weak-bound}
Let $\cT$ be a breadth-critical tangle of order
$q\ge3$ in a matroid $M$, and let $U$ be a
witness for $\br_M(\cT)$.
Every $\cT$-weak set
$H\subseteq E(M)\setminus U$ satisfies
\[
|H|\le a_{\lambda_M(H)},
\]
where the sequence $(a_t)_{t\ge0}$ is defined
in~\eqref{eq:recurrence}.
\end{proposition}

\begin{proof}
First suppose that $H\subseteq E(M)\setminus U$ is
$\cT$-weak and titanic, and put $t=\lambda_M(H)$.
By Proposition~\ref{prop:compression}, there are disjoint
sets $C,D\subseteq H$ such that, for $N=(M/C)\setminus D$
and $H'=H\setminus(C\cup D)$, the set $H'$ is titanic in $N$,
$\lambda_N(H')=t$, and $|H'|\le h_t$.
Since $H$ is $\cT$-weak and $H\cap U=\varnothing$,
Theorem~\ref{thm:transfer} implies that $\cT$ generates a
tangle $\cS$ of order $q$ in $N$ with
$\br_N(\cS)=\br_M(\cT)$.
If $|H|>h_t$, then $|H'|<|H|$, so $C\cup D\ne\varnothing$
and $N$ is a proper minor of $M$.
This contradicts breadth-criticality of $\cT$.
Therefore $|H|\le h_t$ whenever $H$ is titanic.

We now prove the proposition by induction on $t=\lambda_M(H)$,
with $M$, $\cT$, and $U$ fixed.
If $t=0$, then $H$ is titanic because connectivity is nonnegative.
The preceding argument gives $|H|\le h_0=a_0=0$.

Suppose that $t\ge1$ and that the proposition holds for every
$\cT$-weak set $J\subseteq E(M)\setminus U$ with
$\lambda_M(J)<t$. Let $H\subseteq E(M)\setminus U$ be
$\cT$-weak with $\lambda_M(H)=t$.
If $H$ is titanic, the bound already proved gives
$|H|\le h_t\le a_t$.
Otherwise, there is a partition $H=H_1\sqcup H_2\sqcup H_3$
with $\lambda_M(H_j)<t$ for every $j\in\{1,2,3\}$.
Since $H$ is $\cT$-weak, there exists
$F\in\cT$ such that $H\subseteq F$.
For every $j\in\{1,2,3\}$, we have $H_j\subseteq F$ and
$H_j\cap U=\varnothing$. Thus every $H_j$ is $\cT$-weak
and satisfies the induction hypothesis.
Since $(a_t)_{t\ge0}$ is nondecreasing and
$a_t=\max\{h_t,3a_{t-1}\}$, we obtain
\[
|H|=\sum_{j=1}^{3}|H_j|
\le\sum_{j=1}^{3}a_{\lambda_M(H_j)}
\le3a_{t-1}\le a_t.
\]
This completes the induction.
\end{proof}

\begin{proof}[Proof of Theorem~\ref{thm:main}]
Fix an integer $i\ge0$, and let $\cT$ be a
breadth-critical tangle of order $q\ge i+2$
in a finite matroid $M$.
Let $A\subseteq E(M)$ satisfy $\lambda_M(A)=i$.

If $i=0$, Lemma~\ref{lem:critical-connectivity}
implies that $M$ is connected.
Since $\lambda_M(A)=0$, we have
$A=\varnothing$ or $A=E(M)$.
Hence
\[
\min\{|A|,|E(M)\setminus A|\}=0=s_0.
\]

Suppose that $i\ge1$. Then $q\ge3$.
Choose a witness $U$ for $\br_M(\cT)$.
Since $i<q-1$, axiom~(T2) gives a set
$F\in\{A,E(M)\setminus A\}$ with $F\in\cT$.
Put $F^c=E(M)\setminus F$.
By symmetry, $\lambda_M(F)=i$.

Put $I=F\cap U$ and $H=F\setminus U$.
Then $F=I\sqcup H$ and
$E(M)\setminus H=F^c\cup I$.
Since $M_{\cT}|U$ is a spanning uniform
restriction and $F\in\cT$, we have
\[
\min\{q-1,|I|\}=\rho_{M,\cT}(I)
\le\rho_{M,\cT}(F)
\le\lambda_M(F)=i<q-1.
\]
Hence $|I|\le i$.
Since $H\subseteq F$, we have
$r_M(H)\le r_M(F)$.
Moreover, submodularity of the rank function gives
$r_M(F^c\cup I)\le r_M(F^c)+r_M(I)
\le r_M(F^c)+|I|$.
Therefore,
\begin{align*}
\lambda_M(H)
&=r_M(H)+r_M(F^c\cup I)-r(M)\\
&\le r_M(F)+r_M(F^c)+|I|-r(M)\\
&=\lambda_M(F)+|I|\\
&=i+|I|\le2i.
\end{align*}
Since $H\subseteq F\in\cT$, the set $H$
is $\cT$-weak. Moreover, $H\cap U=\varnothing$.
By Proposition~\ref{prop:weak-bound} and the
monotonicity of $(a_t)$,
\[
|F|=|I|+|H|
\le i+a_{\lambda_M(H)}
\le i+a_{2i}=s_i.
\]
Thus
\[
\min\{|A|,|E(M)\setminus A|\}\le s_i.
\]

Finally, let $k\ge2$, and let $\cT$ be a
breadth-critical tangle of order $q\ge k$
in a finite matroid $M$.
For every $i\in\{0,1,\ldots,k-2\}$, we have
$q\ge i+2$.
The bound proved above therefore shows that
$M$ is $(s_0,s_1,\ldots,s_{k-2})$-connected.
\end{proof}

The known bounds $0$, $1$, and $4$ can replace the first
three terms of the sequence $(s_i)_{i\ge0}$.

\begin{corollary}\label{cor:initial}
Conjecture~\ref{conj:81} holds with the sequence
$(\widehat s_i)_{i\ge0}$ defined by $\widehat s_0=0$,
$\widehat s_1=1$, $\widehat s_2=4$, and
$\widehat s_i=i+a_{2i}$ for every integer $i\ge3$.
\end{corollary}

\begin{proof}
Fix an integer $k\ge2$, and let $\cT$ be a
breadth-critical tangle of order $q\ge k$
in a finite matroid $M$.
Let $i\in\{0,\ldots,k-2\}$, and let
$A\subseteq E(M)$ satisfy $\lambda_M(A)=i$.

If $i=0$, then $M$ is connected by
Lemma~\ref{lem:critical-connectivity}.
Hence $A=\varnothing$ or $A=E(M)$, so
$\min\{|A|,|E(M)\setminus A|\}=0=\widehat s_0$.

If $i=1$, then $q\ge3$.
By Lemma~\ref{lem:critical-connectivity},
$M$ is $3$-connected. Thus
$\min\{|A|,|E(M)\setminus A|\}\le1=\widehat s_1$.

If $i=2$, then $q\ge4$.
By Lemma~\ref{lem:critical-connectivity},
$M$ is weakly $4$-connected. Hence
$\min\{|A|,|E(M)\setminus A|\}\le4=\widehat s_2$.

Finally, suppose that $i\ge3$.
Since $q\ge k\ge i+2$,
Theorem~\ref{thm:main} gives
\[
\min\{|A|,|E(M)\setminus A|\}
\le s_i=\widehat s_i.
\]
Therefore $M$ is
$(\widehat s_0,\widehat s_1,\ldots,
\widehat s_{k-2})$-connected.
\end{proof}

\end{document}